\documentclass[openacc]{rsproca_new}%
\usepackage{color}
\usepackage{enumitem}
\usepackage{lineno}
\usepackage{amsmath}
\usepackage{bm,upgreek}
\usepackage{cite}
\usepackage{ragged2e}
\usepackage{changes}
\usepackage{lipsum}
\usepackage{pifont}
\usepackage{subfig}
\usepackage{amsthm}

\usepackage{graphics}
\usepackage{graphicx}
\usepackage{mathtools}

\allowdisplaybreaks[4]

\newtheorem*{definition*}{Definition}

\definecolor{DXH}{rgb}{0,0,1}
\definecolor{DXHong}{rgb}{0.7,0,0.2}

\newcommand{\mylabel}[2]{#2\def\@currentlabel{#2}\label{#1}}

\definechangesauthor[name={DXH}, color=orange]{dxh}
\definechangesauthor[name={DXH}, color=magenta]{dxh2}

\begin{document}

\title{Locating isolas in nonlinear oscillator systems using uncertainty quantification}

\author{%
Dongxiao Hong$^1$, David A.W.~Barton$^2$ and Simon A.~Neild$^2$}

\address{$^1$Department of Mechanical Engineering, Imperial College London, SW7 2AZ, UK

	    $^2$Faculty of Science and Engineering, University of Bristol,  BS8 1TR, UK}

\subject{mechanical engineering}


\corres{Dongxiao Hong\\
\email{d.hong@imperial.ac.uk}}

\begin{abstract}
\sloppy
Parametric uncertainty in nonlinear dynamical systems can fundamentally alter bifurcation behaviour, leading to qualitative response changes.
Predicting operating margins/envelopes under such uncertainties is critical but challenging: conventional uncertainty quantification (UQ) methods struggle to efficiently propagate uncertainties across bifurcation boundaries, where response gradients become singular and solution branches emerge/vanish.

We present a general UQ framework for bifurcation analysis of nonlinear dynamical systems with proportional parametric uncertainty, which systematically integrates continuation methods with parametric sensitivities and extremal conditions.
The approach uses a two-step scheme: first, the loci of extremal response points are traced as the uncertainty domain is expanded from a deterministic reference point; second, these extremal points are tracked as the bifurcation parameter varies, thus determining the maximum and minimum response margins throughout.
The continuation problem scales linearly with the number of uncertain parameters, enabling efficient analysis.

The method is demonstrated on a two-degree-of-freedom nonlinear oscillator exhibiting a range of bifurcation phenomena, including multiple solutions, modal interactions, and symmetry breaking.
In all cases, the framework efficiently captures uncertainty-induced shifts in bifurcation boundaries and response margins.
Notably, the method reveals that parametric uncertainty induces topological changes in the bifurcation structure, including the emergence of an isolated response branch that is absent in the deterministic system with the reference parameters.

\end{abstract}
\begin{fmtext}
\end{fmtext}
\maketitle
\newpage
\noindent
\section{Introduction}\label{sec:intro}

Parametric uncertainty can significantly affect dynamical system behaviour, potentially leading to unanticipated regimes that cause irreversible damage or catastrophic failure.
A notable engineering example is the self-sustaining, high-frequency shimmy in landing gear, triggered by uncertain tyre wear, environmental conditions, or dynamic loads \cite{Tartaruga15,Yin16}.
To mitigate the risk of damage or failure, Uncertainty Quantification (UQ) is essential.

A substantial body of UQ methodologies has been developed and applied across various fields \cite{Marelli14,Sudret17,Chao23}.
A common UQ approach is the Monte Carlo (MC) method \cite{Nicholas49}, which relies on random realisations of the uncertainty.
Despite its conceptual simplicity, MC typically requires a large number of samples for convergence.
To mitigate this cost, stochastic collocation (SC), and polynomial chaos expansion (PCE) can be employed;
non-intrusive SC and PCE employ prescribed samples to construct an interpolant or an orthogonal polynomial expansion for uncertainty propagation \cite{Wiener38,Xiu02,Loeven07,Michael09a}.
When probabilistic information is limited or unnecessary, interval analysis offers a practical alternative \cite{Moore79,Worden05,Chao20}.
Sensitivity-based analysis leverages local information to assess the impact of uncertainty \cite{Tortorelli94,Wilkins09,Dankowicz21}, and is widely used in uncertainty-aware optimisation \cite{Rubino18,Li20}.

Whilst these UQ methodologies have proven effective for linear systems, extending them to nonlinear regimes is nontrivial because bifurcations, i.e.,~qualitative or topological changes in dynamics induced by perturbations, can fundamentally alter the response characteristics \cite{Nayfeh89,Rodrigues10,Botts12,Kuether15,Detroux15,Mingwu22,Dongxiao22a,Lamb26}.
Given that bifurcations can correspond to critical, potentially catastrophic events, quantifying their onset in the presence of uncertainty is crucial.
The probability of the occurrence and type of bifurcations has been investigated by combining centre manifold reduction and coefficient analysis in \cite{Kuehn21}.
The minimum uncertain perturbation that makes bifurcation occur can be determined by constructing an optimisation problem \cite{Iannelli21}.
From a set-valued perspective on uncertain dynamical systems, the bifurcation scenarios have been mathematically distinguished via global analysis in \cite{Lamb15, Lamb26}.

In addition to detecting the existence of bifurcations, locating bifurcation margins is particularly valuable because these margins often define the admissible operating envelopes for uncertain systems.
Approaches based on MC, SC, and PCE have been successfully employed for this task \cite{Attar05,Philip06,Tartaruga18,Breden20,Kuehn24}.
Nonetheless, substantial implementation challenges arise in the presence of intricate bifurcation phenomena such as multiple solutions.
In such cases, constructing a global single-valued interpolant or polynomial expansion for uncertainty propagation is difficult and often infeasible \cite{Worden05,Vineeth13,Breden20}. 
Inappropriate implementation may lead to spurious margins and unreliable predictions \cite{Sinou152,Panunzio17}.
To handle this situation, a multi-resolution PCE approach can be employed by introducing new variables and subdividing the uncertainty domains \cite{Wan06}.
Another PCE-based approach considers intrusive frameworks combining PCE with harmonic balance \cite{Didier13,Sinou15}, treating the bifurcation parameter as stochastic.
Alternatively, special projections can locally restore a single-valued mapping by introducing a fixed parameter, e.g., arc-length ratio \cite{Panunzio17}, solution interval \cite{Tartaruga18}, polar angle \cite{Chao22}, or response phase \cite{Roncen18}.
Despite successes, the lack of a general rule for selecting a suitable projection necessitates a careful examination of the bifurcation diagram before conducting UQ, restricting these approaches to case-by-case studies.
Alternatively, a boundary map can be employed, which corresponds in a unique way to the boundaries of attractors of the multi-solution dynamical system \cite{Lamb23, Kourliouros23, Lamb26}.
Besides implementation challenges, most approaches rely on resolving the output surface over the \emph{entire} uncertainty domain to determine output margins. 
Consequently, the computational cost increases rapidly with the number of uncertain parameters, leading to the well-known curse of dimensionality \cite{Xiu02, Panunzio17}.

We address these challenges by developing a UQ framework for nonlinear dynamical systems with bounded parametric uncertainties.
In contrast to conventional UQ workflows that first reconstruct the output surface and subsequently extract margins, our approach formulates the UQ task directly as a margin-characterisation problem.
The output margins and the associated margin-governing parametric uncertainties are obtained by solving a constrained extremum problem.
Importantly, this formulation avoids the high computational cost of resolving high-dimensional output surfaces;
instead, it reduces the UQ task to tracking one-dimensional margins.
Consequently, the computational cost scales linearly with the number of uncertain parameters, enabling efficient UQ analysis.
Specifically, the problem is solved via numerical continuation, of which the idea is tracking solutions based on local sensitivity by a predictor-corrector scheme \cite{Allgower12, Dankowicz13, Detroux15, Grenat19}.
It should be noted that continuation methods are extensively used in UQ studies, primarily as a computational tool for the realisations of sample-based SC and PCE methods \cite{Didier13,Sinou15,Sinou152,Panunzio17,Tartaruga18,Chao22,Roncen18}.
Alternatively, continuation has been implemented to track local statistical quantities of equilibria by deriving an augmented problem for systems with uncertain states \cite{Kuehn12}.
This idea has later been extended to account for limit cycles and tori \cite{Ahsan24}.
Nonetheless, there is little work that leverages continuation methods to track margins in the presence of parametric uncertainty.

To demonstrate the proposed UQ framework, intricate bifurcation phenomena such as multiple solutions, modal interactions, and symmetry breaking are considered in the presence of parametric uncertainty.
Particularly, a challenging nonlinear feature -- the isolated response curve, or \emph{isola}, which is detached from the primary response branch \cite{Alexander09,Gatti11,Kuether15,Renson19, Donmez25}, is investigated.
Detecting isolas typically requires specialist knowledge and substantial computational effort due to their detached features.
Primary approaches identify a seed point via bifurcation theory, then use numerical continuation to uncover the full branch \cite{Kuether15,Hill16,Grenat19,Dongxiao22, Donmez25}.
This continuation methodology was later applied in the experimental context, termed control-based continuation, for tracking isolas \cite{Renson19}.
As complementary information, singularity theory has been used to provide topological characterisation of isolas \cite{Habib18}.
Alternatively, global analysis has been used to detect an isola by finding initial conditions within the basin of attraction of that isola \cite{Noel15}.
However, to the best of the authors’ knowledge, locating isolas in the presence of uncertainty remains an open question.

The rest of the paper is organised as follows.
$\S\ref{sec:moti}$ introduces the concept via a motivating example.
$\S\ref{sec:framework}$ derives the mathematical formulation, comprising equational constraints from the dynamics, extremal conditions, and parametric sensitivity; the resulting system is solved using a two-step continuation scheme.
$\S\ref{sec:egs}$ presents four applications of varying complexity, showcasing the effectiveness of the UQ framework in the presence of intricate bifurcation phenomena.

\section{Motivating Example}\label{sec:moti}

We use a Duffing oscillator to illustrate typical nonlinear features and the limitations of MC-based UQ, motivating the proposed method. More insightful examples are presented in later sections. The differential equations governing the dynamics are
\begin{equation}\label{eq:EoM_duffing}
m \ddot{q} + c \dot{q} + k q + \alpha q^3 = F\cos\left(\omega t\right)\,,
\end{equation}
where $q \coloneq q(t;q_0)$ is displacement, with velocity $\dot{q}$ and acceleration $\ddot{q}$;
$m$, $c$, $k$, and $\alpha$ are the mass, damping, linear stiffness, and cubic stiffness;
$F$ and $\omega$ are the excitation amplitude and frequency.
The reference parameters (i.e., with no uncertainty) are $m = 1$, $c=0.1$, $k=1$, $\alpha=1$, and $F=0.2$.

\begin{figure}[t!]
   \centering
   \includegraphics[width=9cm]{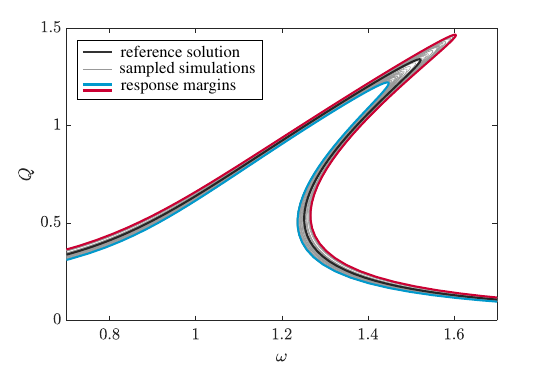}
   \caption{	Uncertainty quantification for the Duffing oscillator.
		The forced response curve of the reference system (without uncertainty) is shown as a black line in the projection of excitation frequency against displacement amplitude.
		The response margins of the uncertain system are obtained via a simulation-based UQ method.
		(Online version in colour.)}
   \label{fig:duffing-ref-simu}
\end{figure}

Under periodic excitation, the system exhibits a family of periodic orbits (POs), commonly visualised as a forced response curve (FRC) in the amplitude–frequency plane;
the reference FRC is the black curve in Fig.~\ref{fig:duffing-ref-simu}.
Nonlinear hardening bends the curve towards higher frequencies, producing multiple solutions near $\omega \approx 1.4$.

In UQ, a common goal is to determine output margins that bound the system’s response under parametric uncertainty. MC achieves this by sampling uncertain parameters from a prescribed domain, computing responses, and inferring extrema. For illustration, we introduce uncertainty in damping and forcing amplitude, compute FRCs for sampled parameter sets (grey curves), and infer response bounds (blue and red) in Fig.~\ref{fig:duffing-ref-simu}.\footnote{The focus here is to provide a straightforward implementation of the UQ approach;
for the sake of simplicity, the detailed specification of uncertain parameters and sampling strategy is not discussed here.}

There are several limitations in applying such a UQ methodology.
Firstly, relying exclusively on simulations to evaluate output margins can be computationally prohibitive, given the large number of required samples.
Apart from determining output margins, assessing the relative importance of uncertain parameters is crucial for understanding the system's performance, which, however, requires complementary studies on uncertainty propagation.
Yet, due to the bifurcation-induced multiple solutions, implementing other approaches, e.g.,~SC and PCE, is challenging and often infeasible \cite{Worden05,Didier13,Sinou152,Panunzio17}.

We propose a UQ framework to determine output margins and margin-governing parametric uncertainties in nonlinear dynamical systems with bounded parametric uncertainty.
The UQ problem is formulated by integrating dynamics, extremal conditions, and parametric sensitivities, and is solved in two steps:
\begin{enumerate}[nosep]
\item Expansion step: trace extremal conditions as the uncertainty domain expands from zero up to a marginal radius, moving the reference solution to a boundary extremum.
\item Propagation step: propagate uncertainty through the bifurcation parameter to trace the margins across the parameter range.
\end{enumerate}
Both steps use numerical continuation as the primary computational tool.

In Fig.~\ref{fig:duffing-ref-simu}, the expansion step starts from a reference solution (a point on the black curve) and advances to the margins (points on the red and blue curves). The propagation step continues these marginal points over frequency, yielding the full margins and associated uncertainties.

\section{Mathematical Framework}\label{sec:framework}

This study considers POs of a general nonlinear dynamical system whose dynamics is governed by a set of first-order ordinary differential equations \cite{Wilkins09,Alexander09,Detroux15,Roncen18,Reddy19,Mingwu22}, i.e.,
\begin{equation}\label{eq:EoM}
\frac{\rm d }{{\rm d} t} \mathbf{x}(t, \mathbf{p}; \mathbf{x}_0 )=
\mathbf{F}\left(\mathbf{x,p}\right)\,,
\end{equation}
where $\mathbf{x} \in \mathbb{R}^{N}$ is the state vector of the system,
$\mathbf{F}$ is a set of nonlinear autonomous functions\footnote{To account for POs of more general periodically-driven  systems, one can augment the underlying autonomous systems with additional state variables to describe the non-autonomous periodic dependency \cite{Dankowicz13}.},
$\mathbf{x}_0\coloneq\mathbf{x}\left(t=0\right)$ is the initial state vector,
$\mathbf{p}$ represents parameters of the system and is described by function
\begin{equation}\label{eq:uncertainty}
\mathbf{p}=\mathbf{P}(\mathbf{p}_0,\mathbf{\epsilon})\,,
\end{equation}
where $\mathbf{p}_0\in\mathbb{R}^{P}$ and $\mathbf{\epsilon}\in\mathbb{R}^{M}$ represent the reference parameter and parametric uncertainty respectively.
A common form of parametric uncertainty is proportional uncertainty, described by a scaling factor of the reference parameter, namely, $\mathbf{p} = \mathbf{p}_0\left(1 + \mathbf{\epsilon} \right)$.

For convenience, the uncertainty is assumed to be bounded by an enclosure
\begin{equation}\label{eq:sphere-constraint}
f_{\epsilon}=\sum_{m} \epsilon_m^2 - r^2 = 0\,,
\end{equation}
where $r\in[0,~R]$ represents the total level of uncertainty.
Eq.~$(\ref{eq:sphere-constraint})$ places each of the uncertain parameters within a hyper-sphere that determines the maximum level of uncertainty for each parameter.
The reference solution, with respect to reference parameters (with no uncertainty), is obtained at $r=0$, whilst the marginal solution, corresponding to the marginal level of uncertainty, is achieved at $r=R$.
In practice, the levels of uncertainty may be informed by expert judgement, experimental testing, and statistical inference from available data.
In other scenarios such as design settings, however, they need not be prescribed a priori: a specification on the desired output margin can instead be used to infer the admissible uncertainty level required to satisfy the performance constraint.
Constraint~$(\ref{eq:sphere-constraint})$ couples the uncertainty in each parameter, restricting the maximal parameter values that can be achieved simultaneously.
This choice is for ease of implementation and alternative constraints are allowable provided that they enclose the domain with a suitably smooth boundary.

The use of a non-smooth boundary function (e.g., a hyper-cube) may also possible.
However, it would likely require extensive modifications to existing numerical continuation codes to deal with transitions across non-smooth points in piecewise-defined constraints.
As such, this is beyond the scope of this paper.

In the presence of uncertainty, the performance of a dynamical system is primarily evaluated based on its uncertain state information.
In practical applications, a scalar-valued performance metric (output) is employed to provide a tractable and application-oriented characterisation.
In general, this metric can be written as
\begin{equation}\label{eq:metric}
g(t,\mathbf{x,p};\mathbf{x}_0):~\mathbb{R}^{N+P+M}\rightarrow\mathbb{R}\,,
\end{equation}
an example of which is the Duffing oscillator's displacement in the motivating example in $\S$\ref{sec:moti}.

For system~$(\ref{eq:EoM})$ equipped with uncertainty~$(\ref{eq:uncertainty})$, defined over domain~$(\ref{eq:sphere-constraint})$, and performance metric~$(\ref{eq:metric})$, the objective of the UQ study is to determine \emph{how uncertainty propagates through the system and leads the reference solution (with no uncertainty) to the extrema on the margin}.

\subsection{Periodic Orbit under Extremal Conditions}\label{sec:constraints}
The periodicity of POs in a nonlinear system~$(\ref{eq:EoM})$ is expressed as an initial-value problem that enforces the recurrence condition over a single period
\begin{equation}\label{eq:BVE}
\mathbf{f}_t = \mathbf{x}\left(T,\mathbf{p}; \mathbf{x}_0\right) - \mathbf{x}_0 = \mathbf{0}\,,
\end{equation}
where $T$ is the response period.%

In addition to periodic conditions~$(\ref{eq:BVE})$, a phase-locking condition is required to fix infinitely many solutions to an isolated point on the PO.
This condition can be arbitrarily defined as long as it imposes a valid phase constraint on states.
Here, we leverage the arbitrariness and define the phase-locking condition as the extremal condition (over a period) on the system's metric.
Particularly, for a system exhibiting a PO, the metric~$(\ref{eq:metric})$ is also a PO whose extremal condition over a period is governed by zero time derivative, i.e.,
\begin{equation}\label{eq:least_favourable_con_time}
L_t\left(0, \mathbf{x,p};\mathbf{x}_0\right)=\frac{{\rm d}g\left(0, \mathbf{x,p};\mathbf{x}_0\right)}{{\rm d}t}=0\,.
\end{equation}
Here, as a common practice, the phase-locking condition is applied to the initial state at $t=0$.

In the presence of uncertainty, apart from the extremal condition~$(\ref{eq:least_favourable_con_time})$ over a period, another extremal condition can be sought in the uncertainty domain.
For an arbitrarily given uncertainty level $r$, the derivative of the metric (with respect to the uncertainty) along the tangent direction of the uncertainty must be zero when extremal uncertainties are achieved.
This is governed by
\begin{equation}\label{eq:least_favourable_con_uncertainty}
\mathbf{L}_{\mathbf{\epsilon}}\left(0, \mathbf{x,p};\mathbf{x}_0\right)=
\Delta_r^{\rm T}
\cdot
\left[
\frac{\partial g}{\partial \mathbf{x}}\frac{\partial \mathbf{x}}{\partial \mathbf{x}_0}\frac{\partial \mathbf{x}_0}{\partial \mathbf{p}}\frac{\partial \mathbf{p}}{\partial \mathbf{\epsilon}}	 +
\frac{\partial g}{\partial \mathbf{p}}\frac{\partial \mathbf{p}}{\partial \mathbf{\epsilon}}
\right] =
\mathbf{0}\,,
\end{equation}
where $\Delta_r^{\rm T}$ is the tangent direction, or null space of the Jacobian, of the uncertainty domain with a level of uncertainty $r$,
and the expression in the brackets is the partial derivative of the metric with respect to uncertainty.
By definition, $\frac{\partial g}{\partial \mathbf{x}}$ and $\frac{\partial g}{\partial \mathbf{p}}$ are analytically defined by the metric expression~$(\ref{eq:metric})$, whilst $\frac{\partial \mathbf{p}}{\partial \mathbf{\epsilon}}$ is analytically defined by the form of uncertainty~$(\ref{eq:uncertainty})$;
$\frac{\partial \mathbf{x}}{\partial \mathbf{x}_0}\left(0, \mathbf{x,p};\mathbf{x}_0\right)$ is an $N \times N$ identity matrix, $\mathbb{I}^{N \times N}$, and $\frac{\partial \mathbf{x}_0}{\partial \mathbf{p}}\frac{\partial \mathbf{p}}{\partial \mathbf{\epsilon}}$ represents the parametric sensitivity of the initial state with respect to uncertainty.

To assess extremal parametric condition~$(\ref{eq:least_favourable_con_uncertainty})$, the parametric sensitivity of the initial state, $\frac{\partial \mathbf{x}_0}{\partial \mathbf{p}}\frac{\partial \mathbf{p}}{\partial \mathbf{\epsilon}}$, is still required.
It can be derived by carrying out sensitivity analysis of dynamics~$(\ref{eq:EoM})$, $(\ref{eq:BVE})$, and $(\ref{eq:least_favourable_con_time})$.
To preserve the clarity and continuity of the discussion, the sensitivity analysis is provided in Appendix~$\ref{sec:appendixA}$.

The above derivations define the POs under extremal conditions through equational constraints~$(\ref{eq:sphere-constraint})$, $(\ref{eq:BVE})$, $(\ref{eq:least_favourable_con_time})$, and $(\ref{eq:least_favourable_con_uncertainty})$.
They can be collected as
\begin{equation}
\mathbf{y}(\mathbf{x}_0,~\mathbf{p}_0,~\mathbf{\epsilon},~r,~T)\coloneqq \left[\mathbf{f}_t;~f_{\epsilon};~L_t;~\mathbf{L}_{\epsilon}\right]=\mathbf{0}\,.
\end{equation}
where $\mathbf{y}:~\mathbb{R}^{N+P+M+2}\rightarrow\mathbb{R}^{N+M+1}$ with variables of the initial states, $\mathbf{x}_0\in\mathbb{R}^{N}$, reference parameters, $\mathbf{p}_0\in\mathbb{R}^{P}$, parametric uncertainty, $\mathbf{\epsilon}\in\mathbb{R}^{M}$, level of uncertainty, $r\in\mathbb{R}$, and response period, $T\in\mathbb{R}$.

In the presence of uncertainty, a primary objective is to identify the extremal metric for a system with fixed reference parameters, i.e.,~$\mathbf{p}_0=\text{const.}$, where the problem can be reduced to $\mathbf{y}(\mathbf{x}_0;~\mathbf{\epsilon};~r;~T) = \mathbf{0}$.
Nonetheless, of greater insight and complexity is the investigation on how uncertainty propagates through a varying free parameter, $\lambda \subseteq \mathbf{p}_0$, to the extremal metric -- how the extremal metric changes as the free parameter varies.
This class of study is known as bifurcation analysis, wherein the free parameter is called the bifurcation parameter \cite{Kernevez87,Allgower12,Detroux15}.
For convenience, a single bifurcation parameter is considered and, therefore, variable $\mathbf{p}_0$ is replaced by $\lambda$ and the problem is reduced to
\begin{equation}\label{eq:full_problem}
\mathbf{y}(\mathbf{x}_0,~\mathbf{\epsilon},~r,~T,~\lambda)\coloneqq \left[\mathbf{f}_t;~f_{\epsilon};~L_t;~\mathbf{L}_{\epsilon}\right]=\mathbf{0}\,.
\end{equation}
where $\mathbf{y}:~\mathbb{R}^{N+M+3}\rightarrow\mathbb{R}^{N+M+1}$.
Under the regular value condition, the solution set of Eq.~$(\ref{eq:full_problem})$ represents a two-dimensional manifold embedded in $\mathbb{R}^{N+M+3}$.

\subsection{Continuation-based Implementation}\label{sec:implementation}
\begin{figure}[t]%
    \centering
    \subfloat{\includegraphics[width=6.3cm]{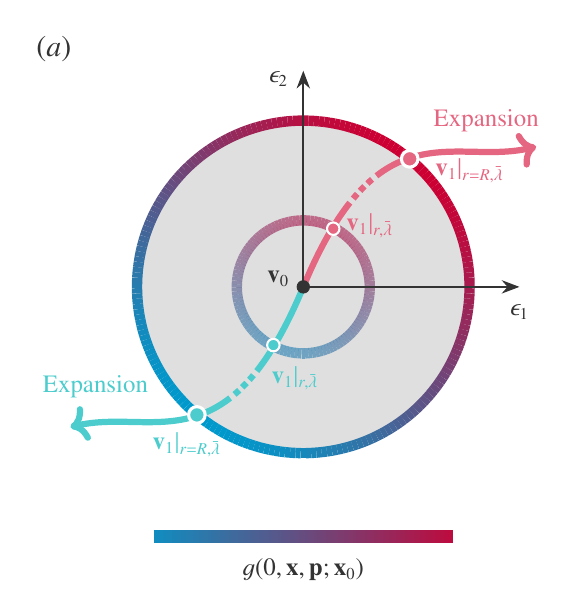} }%
    \quad
    \subfloat{\includegraphics[width=6.7cm]{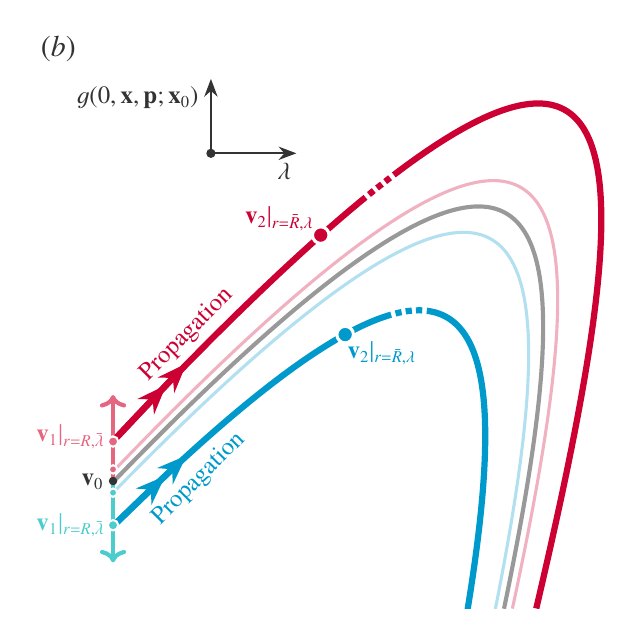} }%
    \caption{	Schematic diagram of the proposed uncertainty quantification method.
		$(a)$~Uncertainty expansion step: the extremal solutions are traced as the uncertainty domain expands from null to the marginal size.
		$(b)$~Uncertainty propagation step: the metric margins are traced as the uncertainty propagates through the bifurcation parameter to the extremal solutions.
		(Online version in colour.)}%
    \label{fig:methodology}%
\end{figure}%

To solve problem~$(\ref{eq:full_problem})$, one can refer to multi-dimensional continuation \cite{Allgower12,Dankowicz13} or successive continuation approaches \cite{Kernevez87,Li18,Li20}.
The following discussions implement the UQ problem via a two-step successive continuation scheme.

\subsubsection*{Step~1: uncertainty expansion}\label{sec:expansion}
As previously mentioned, the proposed method begins with an expansion step, wherein the solution is traced as the uncertainty domain expands from null (with respect to the reference solution with $r=0$) to the marginal level of uncertainty, i.e.,~$r=R$.

This sub-problem is achieved by viewing the level of uncertainty, $r$, as the continuation parameter whilst freezing the bifurcation parameter, $\lambda$, such that problem~$(\ref{eq:full_problem})$ is reduced to %
\begin{equation}\label{eq:propagation_step}
\mathbf{y}_1(\mathbf{x}_0,~\mathbf{\epsilon},~r,~T,~\bar{\lambda})\coloneqq \left[\mathbf{f}_t;~f_{\epsilon};~L_t;~\mathbf{L}_{\epsilon}\right]=\mathbf{0}\,,
\end{equation}
where the frozen variable is denoted by $\bar{\bullet}$ and $\mathbf{y}_1:~\mathbb{R}^{N+M+2}\rightarrow\mathbb{R}^{N+M+1}$ with variables $\mathbf{v}_1 = [\mathbf{x}_0;~\mathbf{\epsilon};~r;~T]$.
This expansion step can be illustrated in the uncertainty domain -- see Fig.~\ref{fig:methodology}$(a)$ for an example case with two uncertain parameters.
The grey region represents the enclosed domain of parametric uncertainty~$(\ref{eq:sphere-constraint})$ where the reference solution is at the origin and denoted $\mathbf{v}_0$ (with reference metric $g_0$), whilst the colour scale is used to denote the metric values.

The expansion step makes use of the analytically derived conditions~$(\ref{eq:propagation_step})$ and traces these conditions from the reference solution as the level of uncertainty increases, shown as one-headed arrows in Fig.~\ref{fig:methodology}.
The solution trajectories are termed positive and negative expansions and denoted in red and blue, respectively. 
For this to be successful, some regular dependence of the continuation problem on the level of uncertainty is required \cite{Allgower12,Dankowicz13}.
Firstly, $\mathbf{y}_1$ is assumed to vary smoothly with $r$ so that it lies on a well-defined local branch and the derivative exists.
It further requires that no bifurcation occurs in the immediate vicinity.
In practice, this does not appear to be a significant restriction.

The objective of this step is consistent with that of UQ methods such as SC and PCE \cite{Wiener38,Xiu02,Loeven07,Michael09a} -- identifying the extremal metric within a prescribed domain of uncertainties.
The primary distinction lies in that SC and PCE approximate the \emph{entire metric surface} over the uncertainty domain using sampling-based functional interpolations;
whereas the proposed step traces only the subset in which conditions~$(\ref{eq:propagation_step})$ hold.

\subsubsection*{Step~2: uncertainty propagation}\label{sec:propagation}
The uncertainty propagation step aims to propagate uncertainty through the bifurcation parameter to solution points from Step~1, as such identifying the margins across the bifurcation parameter and the corresponding margin-governing uncertainties.

This sub-problem is achieved by viewing the bifurcation parameter, $\lambda$, as the continuation parameter and freezing the level of uncertainty to its marginal level, $R$, such that problem~$(\ref{eq:full_problem})$ remains a problem being underdetermined by one degree of freedom,
\begin{equation}\label{eq:expansion_step}
\mathbf{y}_2(\mathbf{x}_0,~\mathbf{\epsilon},~\bar{R},~T,~\lambda)\coloneqq \left[\mathbf{f}_t;~f_{\epsilon};~L_t;~\mathbf{L}_{\epsilon}\right]=\mathbf{0}\,,
\end{equation}
where $\mathbf{y}_2:~\mathbb{R}^{N+M+2}\rightarrow\mathbb{R}^{N+M+1}$ and the variables are $\mathbf{v}_2 = [\mathbf{x}_0;~\mathbf{\epsilon};~T;~\lambda]$.
This propagation step can be illustrated in the projection of the bifurcation parameter against the system metric -- see Fig.~\ref{fig:methodology}$(b)$.
In this projection, the results of the expansion step are firstly revisited --- they are along the metric axis because the bifurcation parameter is fixed.
Initialising from the results of the expansion step, the propagation step is represented by double-headed arrows that continue the metric margins in the bifurcation parameter.

\subsubsection*{Continuation algorithm}\label{sec:continuation}
The above two-step framework yields two zero problems~$(\ref{eq:propagation_step})$ and $(\ref{eq:expansion_step})$, both underdetermined by one degree of freedom, i.e.,~$\mathbf{y}_i:~\mathbb{R}^{N+M+2}\rightarrow\mathbb{R}^{N+M+1}$ where $i = 1,~2$ denotes Step~1 and 2 respectively.
They can be solved via typical one-dimensional numerical continuation algorithms \cite{Allgower12, Dankowicz13, Detroux15, Grenat19}.
The adopted continuation scheme employs a predictor-corrector procedure to trace the solution $\mathbf{v}_i$ for a zero problem $\mathbf{y}_i$.
Here the implementation is briefly discussed for completeness;
for details the reader can refer to \cite{Allgower12, Dankowicz13, Detroux15, Grenat19}.

Starting from a known $j^{\rm th}$ solution\footnote{How to initialise the continuation using a known solution is deferred to the end of this section.}, $\mathbf{v}_{i,j}$, the predictor aims to provide an approximate solution, $\tilde{\mathbf{v}}_{i,j+1}$, along the tangent direction at this known solution.
The tangent vector, $\Delta_{i,j}$, is locally in the null space of the Jacobian at $\mathbf{v}_{i,j}$ and can be solved via
\begin{equation}\label{eq:com-tan}
\left. \dfrac{\partial \mathbf{y}_i}{\partial \mathbf{v}_i} \right|_{\mathbf{v}_i = \mathbf{v}_{i,j}} \cdot \Delta_{i,j} = \mathbf{0}\,.
\end{equation}
Jacobian matrix elements can be computed via finite difference approximations \cite{LeVeque07}.

Equations~$(\ref{eq:com-tan})$ are underdetermined as the number of unknowns outnumbers that of equations.
One typical method to address the problem is fixing one component of the tangent vector and solving the resulting determined equation set.
The obtained tangent vector is then normalised by $\lVert \Delta_{i,j} \rVert = 1$.
Consequently, the predictor, $\tilde{\mathbf{v}}_{i,j+1}$, can be obtained by definition
\begin{equation}
\tilde{\mathbf{v}}_{i,j+1} = \mathbf{v}_{i,j} + h \Delta_{i,j}\,,
\end{equation}
where $h$ denotes the perturbation size along the tangent direction.

Following the predictor is the corrector, which corrects the predictor iteratively along the orthogonal direction to find $\mathbf{v}_{i,j+1}$ on the solution branch.
This orthogonal correction augments the original underdetermined system of equations~$(\ref{eq:propagation_step})$ and $(\ref{eq:expansion_step})$ with one additional constraint and results in a determined problem, i.e.
\begin{equation}
\mathbf{ \mathcal{Y} }_i =
\begin{bmatrix}
\mathbf{y}_i	( \mathbf{v}_{i,j+1} )	\\
\Delta_{i,j}^{\rm T} \cdot \left( \mathbf{v}_{i,j+1} - \tilde{\mathbf{v}}_{i,j+1} \right)
\end{bmatrix}
= \mathbf{0}\,,
\end{equation}
where $\mathbf{ \mathcal{Y} }_i:~\mathbb{R}^{N+M+2}\rightarrow\mathbb{R}^{N+M+2}$ is the augmented $i^{\rm th}$ step problem and the last equation denotes the orthogonal correction.
Using the predictor as the initial solution guess, this system of equations can be solved via iterative methods such as Newton's method.

As discussed earlier, to initialise the continuation, a known solution that satisfies each problem is needed.
For the uncertainty expansion step, this solution is obtained by introducing uncertainty perturbations to the reference system.
This results in a determined problem that can be solved via Newton's method using the reference solution, $\mathbf{v}_0$, as the initial guess.
The results of the uncertainty expansion step, $\mathbf{v}_1\vert_{r=R}$, is then passed to the uncertainty propagation step to continue the margins.

\section{Case Studies of Uncertainty Quantification}\label{sec:egs}

This section presents applications of the developed UQ method.
Here, a two-mode system with cubic nonlinearity is considered, the dynamics of which is governed by a set of non-autonomous equations of motion
\begin{subequations}\label{eq:two-mode-EoM}
\begin{align}
m_1 \ddot{q}_1 					+
(c_1 + c_2) \dot{q}_1				-
c_2 \dot{q}_2					+
(k_1 + k_2) q_1					-
k_2 q_2						+
\alpha_1 q_1^3 + \alpha_2(q_1 - q_2)^3 	=
F_1\cos\left(\omega t \right)\,,				\\
m_2 \ddot{q}_2 					-
c_2 \dot{q}_1					+
(c_2 + c_3) \dot{q}_2				-
k_2 q_1						+
(k_2 + k_3) q_2					+
\alpha_3 q_2^3 + \alpha_2(q_2 - q_1)^3 	=
F_2\cos\left(\omega t \right)\,,
\end{align}
\end{subequations}
where $m_1$ and $m_2$ are mass coefficients, $k_1$, $k_2$, and $k_3$ are linear stiffness coefficients, $c_1$, $c_2$, and $c_3$ are damping coefficients, $\alpha_1$ and $\alpha_2$ are cubic nonlinear stiffness coefficients, $F_1$ and $F_2$ are excitation amplitudes, and $\omega$ is the excitation frequency.

This two-mode model is extensively considered in bifurcation studies of various mechanical systems such as nonlinear cables, beams, and energy sinks \cite{Nayfeh89,Manevitch05,Habib18,Ding20,Dongxiao22}.
Non-autonomous dynamics~$(\ref{eq:two-mode-EoM})$ can be rewritten in an autonomous state-space form as Eqs.~$(\ref{eq:EoM})$ using an augmented state-space description, namely, $\mathbf{x} = [q_1,~\dot{q}_1,~q_2,~\dot{q}_2,~s_1,~s_2]^{\rm T}$, where $s_1$ and $s_2$ represent variables of a Stuart–Landau oscillator (Hopf normal form) and they are used to make the system autonomous by generating the periodic forcing internally \cite{Dankowicz13}.
For details of the rewritten equations, the reader is directed to Appendix~$\ref{sec:appendixB}$.

In the following, four examples are considered, representing uncertain nonlinear systems with different levels of complexity:
\begin{enumerate}[nosep]
\item[$\S 4\ref{sec:duffing-like}$] The system has well-separated eigenfrequencies.
\item[$\S 4\ref{sec:mode-int}$] The system has commensurate eigenfrequencies.
\item[$\S 4\ref{sec:sym}$] The system exhibits symmetry changes in its configurations due to uncertainty.
\item[$\S 4\ref{sec:isola}$] The system exhibits an isolated response curve due to uncertainty.
\end{enumerate}
The reference parameters of the system are listed in Table~$\ref{table:par_eg_sec4}$ whilst different combinations of uncertainty are specified in the subsequent discussions.
For all examples, the forcing frequency, $\omega$, is considered as the bifurcation parameter to investigate the uncertain impact on FRCs.
Here, we note that intricate bifurcations exist for these examples --
they all show multiple solutions, and some present topological changes of FRCs in the presence of uncertainty.
As discussed in $\S\ref{sec:intro}$, these features complicate the implementation of typical UQ methods such as SC and PCE.
Inappropriate implementation may yield spurious margins and make unreliable predictions \cite{Sinou152,Panunzio17}.
Although a variety of special projections have been proposed to address this issue \cite{Panunzio17,Tartaruga18,Roncen18,Chao22}, their applicability is usually case-dependent.
Therefore, instead of relying on typical UQ methods with case-by-case implementations for validation,
all examples are validated via simulations conducted over the gridded-point uncertainty domain.
The corresponding results are provided in Appendix~$\ref{sec:appendixC}$.

\begin{table}[t]
\caption{Parameters of the example systems.}
\label{table:par_eg_sec4}
\centering
\scriptsize
\begin{tabular}{cccccccccccccccc}
	\hline\noalign{\smallskip}
	$\rm Example$			& $m_1$ & $m_2$ & $c_1$ & $c_2$ & $c_3$ & $k_1$ & $k_2$ & $k_3$ & $\alpha_1$ & $\alpha_2$ & $\alpha_3$ & $F_1$ & $F_2$
\\	\noalign{\smallskip}\hline\noalign{\smallskip}
	$\S 4\ref{sec:duffing-like}$	& 1         & 1	        &	 0.05	 & 0.005  & 0.05   &	 1 	& 1 	    & 1        & 1                 & 0.5              &  1                & 0.03    & 0.03
\\	\noalign{\smallskip}\hline\noalign{\smallskip}
	$\S 4\ref{sec:mode-int}$	& 1         & 1	        &	 0.08	 & 0.02   & 0.05    &	 1 	& 0.2 	    & 0.8     & 1                 & 0.02            &  2                & 0.1      & -0.03
\\	\noalign{\smallskip}\hline\noalign{\smallskip}
	$\S 4\ref{sec:sym}$	& 1         & 1	        &	 0.008	 & 0.001   & 0.008    &	 1 	& 0.04 	    & 1     & 0.5                 & 0.01            & 0.5                & -0.005      & 0.0052
\\	\noalign{\smallskip}\hline\noalign{\smallskip}
	$\S 4\ref{sec:isola}$	& 1         & 0.05	        &	 0.015	 & 0.015   & 0       &	 1 	& 0.0454 	    & 0     & 1                 & 0.0042            &  0               & 0.2      & 0
\\        \hline\noalign{\smallskip}
\\
\end{tabular}
\end{table}

\subsection{Well-separated Eigenfrequencies}\label{sec:duffing-like}

\begin{figure}[t]
   \centering
   \includegraphics[width=14cm]{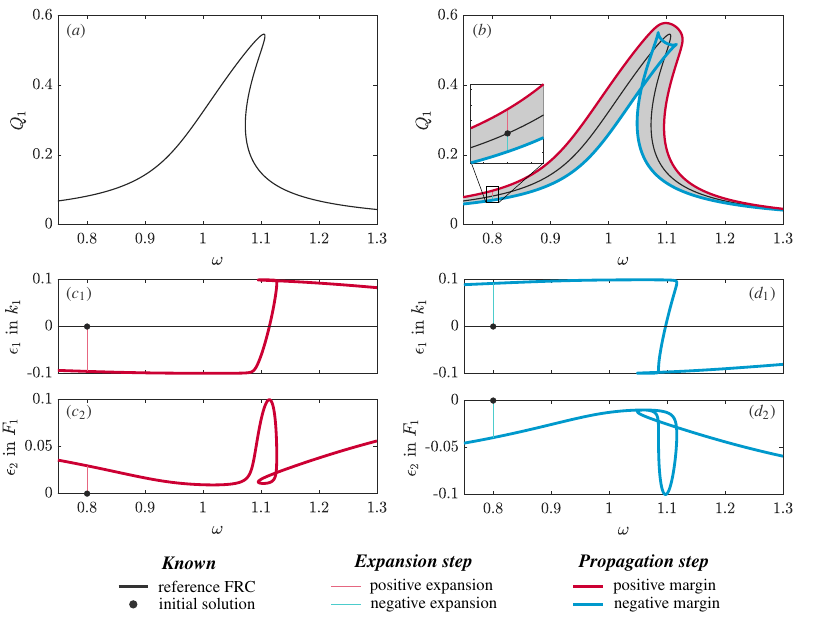}
   \caption{	Uncertainty quantification for a two-mode system with sufficiently separated eigenfrequencies.
		$(a)$~Forced response curve for the reference system (without uncertainty) in the projection of excitation frequency against displacement amplitude of the first mass.
		$(b)$, $(c)$, and $(d)$~Metric margins and margin-governing uncertainties for the uncertain system with an uncertainty level of $R = 0.1$.
		(Online version in colour.)}
   \label{fig:duffing-like}
\end{figure}

The reference parameters of the system are listed in the first row of Table~\ref{table:par_eg_sec4}.
This two-mode system has sufficiently separated eigenvalues ($\omega_{n1} = 1$ and $\omega_{n2}\approx 1.732$), leading to limited interaction between its modal components under excitation.

For this system, the displacement of the first mass is seen as the metric;
and the reference FRC (without uncertainty) is shown in Fig.~\ref{fig:duffing-like}$(a)$ in the projection of excitation frequency against the initial displacement (amplitude) of the first mass.

In the uncertain case, proportional uncertainty is introduced to the linear stiffness and forcing amplitude, specifically, $k_{1,\epsilon} = k_1(1 + \epsilon_1)$ and $F_{1,\epsilon} = F_1(1 + \epsilon_2)$.
The uncertainty is constrained by Eq.~$(\ref{eq:sphere-constraint})$ with a marginal level of $R = 0.1$.
The proposed method is implemented herein to determine the metric margins and their associated extremal uncertainties.
The results are shown in Fig.~\ref{fig:duffing-like}$(a)$, $(b)$, and $(c)$, where the initial solution (with respect to the reference system) is denoted via a black dot, the expansion step results are thin red and blue lines, whilst the propagation step results are thick red and blue lines.

The solutions of the propagation step can be projected onto the excitation frequency against the metric amplitude, see panel~$(b)$, where they represent the margins that bound the metric values in the presence of uncertainty.
It should be noted that, near resonance, the negative margin (thick blue line) manifests a swallowtail catastrophe, arising from including the linear stiffness uncertainty.

Likewise, the solutions can be projected onto the excitation frequency against the uncertainty -- panels~$(c)$ and $(d)$, where they represent the margin-governing uncertainties.
Using these plots, one can compare the relative significance of uncertain parameters.
Away from resonance, the linear stiffness uncertainty shows a dominant impact on the metric, outweighing the effect of the forcing amplitude uncertainty.
In contrast, as the response approaches resonance, the extremal uncertainty undergoes a drastic shift, with dominance transitioning from linear stiffness to forcing amplitude.

\subsection{Commensurate Eigenfrequencies}\label{sec:mode-int}
\begin{figure}[t]
   \centering
   \includegraphics[width=14cm]{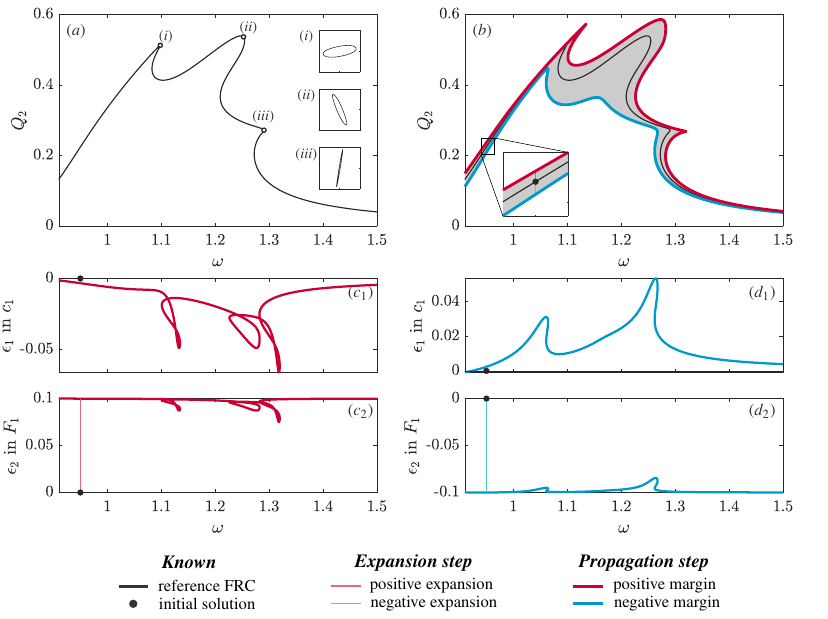}
   \caption{Uncertainty quantification for a two-mode system featuring modal interactions.
		$(a)$~Forced response curve for the reference system in the projection of excitation frequency against displacement amplitude of the second mass.
		The embedded plots~$(i)$, $(ii)$, and $(iii)$ show the periodic orbits for near-resonant POs in the configuration space.
		$(b)$, $(c)$, and $(d)$~Metric margins and margin-governing uncertainties for the uncertain system with an uncertainty level of $R = 0.1$.
		(Online version in colour.)}
   \label{fig:mode-int1}
\end{figure}

\begin{figure}[t]
   \centering
   \includegraphics[width=9cm]{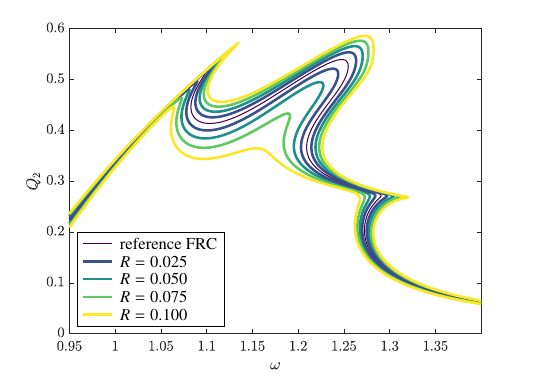}
   \caption{	Sensitivity of the internally resonant region to the level of uncertainty.
		Metric margins for the uncertain system in the projection of excitation frequency against displacement amplitude of the second mass.
		Four uncertainty levels of $R = 0.025,~0.050,~0.075,~0.100$ are shown as thick lines;
		whilst the reference metric is shown as a thin line.
		(Online version in colour.)}
   \label{fig:mode-int2}
\end{figure}

The reference parameters of the second example are given in the second row of Table~\ref{table:par_eg_sec4}.
With these parameters, the system has commensurate natural frequencies ($\omega_{n1} \approx 0.936$ and $\omega_{n2} \approx 1.150$), which can lead to strong modal interactions.
The FRC of the reference system is shown in Fig.~\ref{fig:mode-int1}$(a)$, where the embedded plots display some example POs in the modal configuration space --
the two modes oscillate in a one-to-one frequency commensurability, known as $1:1$ internal resonance \cite{Nayfeh89}.

In this case, uncertain damping coefficient and forcing amplitude are considered, specifically, $c_{1,\epsilon} = c_1(1 + \epsilon_1)$ and $F_{1,\epsilon} = F_1(1 + \epsilon_2)$; the displacement of the second mass, i.e.,~$q_2$, is seen as the metric of interest.
Figure~\ref{fig:mode-int1}$(b)$ presents the metric margins whilst Figs.~\ref{fig:mode-int1}$(c)$ and $(d)$ show the margin-governing uncertainties with an uncertainty level of $R=0.1$.
Panels~$(c)$ and $(d)$ show that, overall, uncertainty in the forcing amplitude exhibits a dominant influence across the frequency range, with small yet sophisticated variations occurring in the vicinity of resonances.

As shown in panel~\ref{fig:mode-int1}$(b)$, the second resonant region has a more extended response range under uncertainty, which indicates a higher sensitivity to uncertainties.
This becomes clearer by considering different levels of uncertainty.
To demonstrate this, the results of four levels of uncertainty, i.e.,~$R = 0.025,~0.050,~0.075,~0.100$ are shown in Fig.~$\ref{fig:mode-int2}$.
The density of margins can serve as a measure of the metric's sensitivity to the level of uncertainty -- the more widely spaced the margins, the greater the sensitivity.

\subsection{Symmetry Changes of Configurations}\label{sec:sym}
\begin{figure}[t]
   \centering
   \includegraphics[width=14cm]{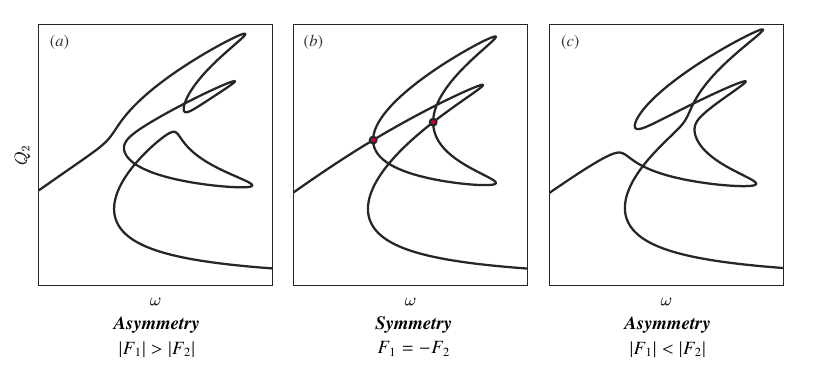}
   \caption{	FRC topologies in the projection of excitation frequency against displacement amplitude of the second mass.
		$(a)$~FRC topology for asymmetric configurations with $\vert F_1 \vert > \vert F_2 \vert$.
		$(b)$~FRC topology for symmetric configurations with $ F_1 =  - F_2 $;
			the topology is featured by two bifurcation points (solid dots) that split when the symmetry breaks.
		$(c)$~FRC topology for asymmetric configurations with $\vert F_1 \vert < \vert F_2 \vert$. }
   \label{fig:FRC_topology}
\end{figure}

\begin{figure}[t]
   \centering
   \includegraphics[width=14cm]{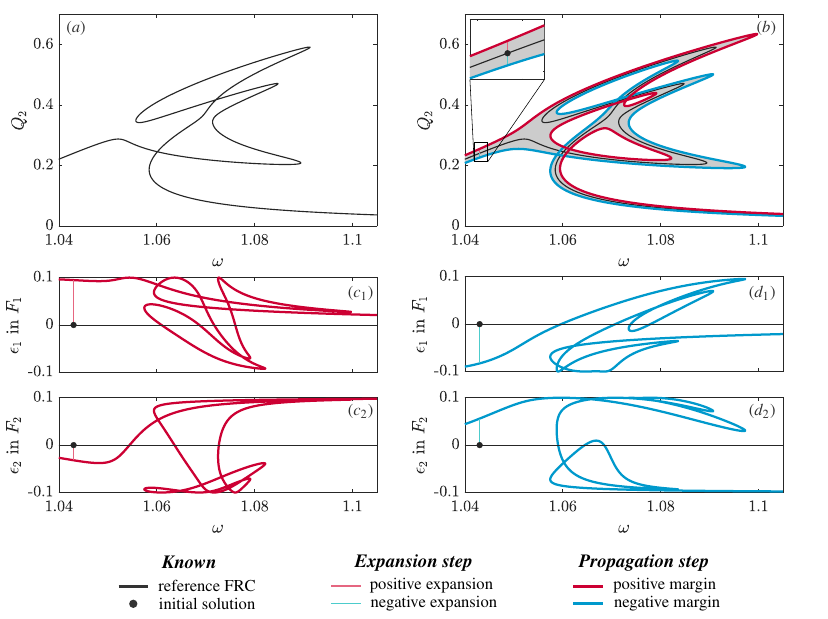}
   \caption{Uncertainty quantification for a system with symmetry changes in its configurations.
		$(a)$~Forced response curve for the reference system in the projection of excitation frequency against displacement amplitude of the second mass.
		$(b)$, $(c)$, and $(d)$~Metric margins and margin-governing uncertainties for the uncertain system with an uncertainty level of $R = 0.1$.
		(Online version in colour.)}
   \label{fig:sym_breaking}
\end{figure}

The reference parameters of the third example are listed in the third row of Table~\ref{table:par_eg_sec4}.
This reference system has an asymmetric configuration, consisting of a symmetric layout, i.e.~$m_1=m_2$, $k_1=k_3$, $c_1=c_3$, and $\alpha_1 = \alpha_3$, and subjected to asymmetric excitations, $F_1\neq -F_2$.

Specifically, the displacement of the second mass is seen as the metric;
additionally, uncertainty is introduced to the asymmetric forcing amplitudes in the form of $F_{1,\epsilon} = F_1(1+\epsilon_1)$ and $F_{2,\epsilon} = F_2(1+\epsilon_2)$ with a marginal level of $R = 0.1$.

Taking the uncertainty into account, the system can show both symmetric ($F_1 = -F_2$) and asymmetric ($F_1\neq -F_2$) configurations.
Correspondingly, there are three different topologies of FRCs for this uncertain system -- see Fig.~\ref{fig:FRC_topology} in the projection of excitation frequency against amplitude of the second mass.
One is related to the symmetric case, featured by two symmetry-breaking bifurcation points -- they are denoted by dots in panel~$(b)$.
As the term indicates, symmetry breaking splits the bifurcation points and leads to FRC topology changes.
There are two FRC topologies for the asymmetric case, depending on the relative quantity of the asymmetric forcing amplitudes.
In cases of $\vert F_1 \vert > \vert F_2 \vert$, the topology is shown in panel~$(a)$, whilst the topology in cases of $\vert F_1 \vert < \vert F_2 \vert$ is shown in panel~$(c)$.
The presence of multiple curve topologies poses significant challenges for implementing UQ methods such as SC and PCE, since interpolation is infeasible in regions with bifurcation splitting.

Here, we demonstrate the applicability of the proposed method in this intricate uncertain scenario.
The system's reference FRC is shown in Fig.~\ref{fig:sym_breaking}$(a)$.
Using the proposed method, extremal solutions are obtained and they are shown as the metric margins in panel~$(b)$.
In this projection, the multi-topology FRC characteristics are captured by margins with different topologies -- the negative margin (blue line) has the same topology as Fig.~\ref{fig:sym_breaking}$(c)$ whilst the positive margin (red line) has the same topology as Fig.~\ref{fig:sym_breaking}$(a)$.
The extremal solutions can alternatively be interpreted as the margin-governing uncertainties, see panels~$(c)$ and $(d)$.
The intricate curves reveal the complexity of extremal uncertain conditions across the frequency range.
Such complexity can pose challenges for sampling-based UQ methods, as resolving the fine structure of the margins may require a prohibitively large number of samples.

\subsection{Locating Isolas}\label{sec:isola}

\begin{figure}[t]
   \centering
   \includegraphics[width=14cm]{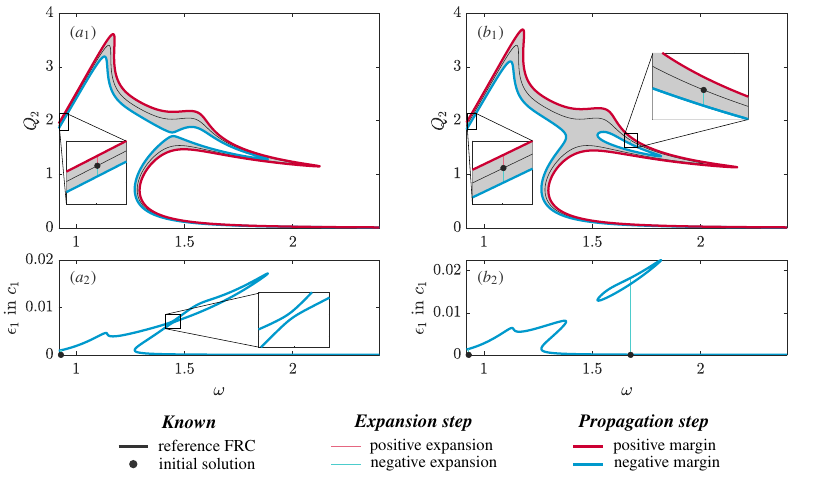}
   \caption{	Locating the isola using uncertainty quantification.
		Metric margins and margin-governing uncertainties for the uncertain system with uncertainty levels of $(a)~R = 0.07$ and $(b)~R = 0.10$.
		(Online version in colour.)}
   \label{fig:isola}
\end{figure}

Nonlinear systems can show rich periodic behaviours, among which the \emph{isola}, i.e.,~an isolated FRC, is of particular importance \cite{Alexander09,Gatti11,Kuether15,Renson19}.
Determining the existence of an isola is generally challenging owing to its being detached from the primary branch.
Much research effort has been dedicated to this problem by seeking a special solution on the isola, following which a continuation algorithm is employed to uncover the full solution set \cite{Kuether15,Noel15,Hill16,Habib18,Grenat19,Dongxiao22}.
To complicate matters, the presence of uncertainty can add to the difficulties in locating isolas.
In this section, the proposed method is applied to handling this task.

For demonstration, the two-mode system~$(\ref{eq:two-mode-EoM})$ is again considered, the reference parameters of which are listed in the fourth row of Table~\ref{table:par_eg_sec4}.
Particularly, the displacement of the second mass is viewed as the metric.
The FRC of the reference system is shown as a black line in Fig.~\ref{fig:isola}, where one single primary branch can be seen.

Uncertain damping coefficient and forcing amplitude are then considered, i.e.,~$c_{1,\epsilon} = c_1(1 + \epsilon_1)$ and $F_{1,\epsilon} = F_1(1 + \epsilon_2)$.
Firstly, a marginal uncertainty level of $R=0.07$ is considered, where the related metric margins and margin-governing uncertainties are shown in Fig.~\ref{fig:isola}$(a)$.
The margins show the same topology as the reference FRC, indicating all uncertain cases share the same FRC topology -- having one single primary branch.
Particularly, the negative margin, i.e.~blue line in panel~$(a_1)$, shows a near self-intersection, where the margin nearly meets but then repels;
similar geometry can be found on the margin-governing uncertainty in panel~$(a_2)$.

The level of the uncertainty is then increased from $R=0.07$ to $R=0.10$ with results presented in Fig.~\ref{fig:isola}$(b)$.
In panel~$(b_1)$, the positive margin (red line) remains similar to that in panel~$(a)$;
however, an isolated negative margin (blue loop) is identified in additional to the primary one (blue line).
This topological change in margins indicates the existence of another FRC topology, namely, having an isola, for the uncertain system.
Likewise, an isolated loop can be found in the margin-governing uncertainties in panel~$(b_2)$.

\begin{figure}[t]
   \centering
   \includegraphics[width=14cm]{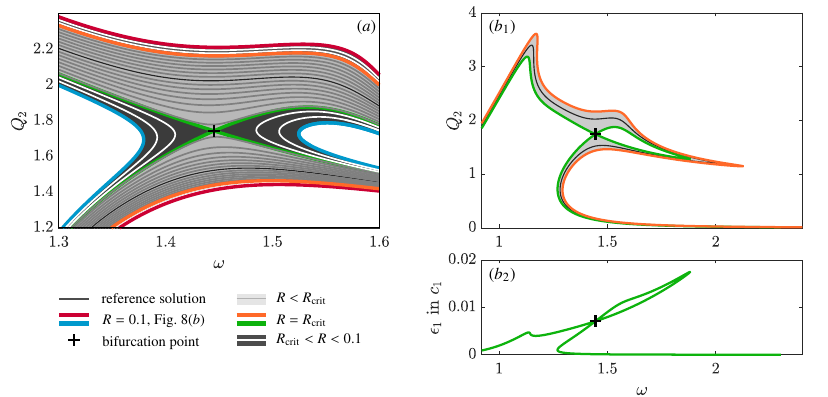}
   \caption{	Tracing the emergence of an isola for an uncertain system.
		$(a)$~Metric margins of cases with different uncertainty levels.
		$(b)$~Critical case with an uncertainty level of $R_{\rm crit}=0.07145$.
		The critical case serves as the condition for the emergence of an isola.
		(Online version in colour.)}
   \label{fig:isola-formulation}
\end{figure}

To trace the topological transition from panel~$(a)$ to panel~$(b)$, cases with different uncertainty levels are shown in Fig.~\ref{fig:isola-formulation}$(a)$ where the margins of the critical case (uncertainty level of $R_{\rm crit} = 0.07145$) are shown as green and orange lines, cases with $R<R_{\rm crit}$ are denoted via grey lines, whilst cases with $R>R_{\rm crit}$ are white lines.
Particularly, the full margins and margin-governing uncertainties of the critical case are shown in panel~$(b)$.
This critical case is characterised by a bifurcation (indicated via a cross) and defines the conditions for the emergence of an isola -- cases with an uncertainty level exceeding $R_{\rm crit}$ may give rise to an isola.

\section{Conclusion}\label{sec:conclusion}
This work presented an uncertainty quantification (UQ) framework for the bifurcation analysis of nonlinear dynamical systems with proportional parametric uncertainty.
The framework is formulated by integrating system dynamics, extremal conditions, and parametric sensitivities.
A two-step successive continuation scheme has been developed to solve the problem and determine the bifurcation diagram margins in the presence of parametric uncertainties.

We demonstrated the method across a range of bifurcation scenarios using a mechanical oscillator.
The first two examples, with well-separated and commensurate eigenfrequencies, show the method’s capability in settings with multiple solutions.
The third example involved uncertainty in system configuration, where the framework identified response margins in the presence of symmetry-breaking bifurcations.
The final example employed the method to locate isolated response curves (isolas) under uncertainty, providing insight into their emergence and evolution.
All examples were validated against simulations, showing excellent agreement.
These results indicate that the framework robustly captures the impact of parametric uncertainty on nonlinear bifurcation behaviour.

We conclude by providing a brief discussion on some potential further research directions connected to this work.
\begin{enumerate}
\item Bifurcations with singularity.
Periodic orbits in nonlinear dynamical systems can approach singular regimes in which standard regular assumptions break down. 
A typical example is near-global bifurcation (e.g.,~homoclinic/heteroclinic connections) where the period grows rapidly and may diverge.
In such situations, the solution branch can become locally non-regular, and the Jacobian becomes ill-conditioned or rank-deficient.
Consequently, continuation-based schemes may fail to provide a reliable analysis near these singular behaviours.
Nonetheless, the proposed method may still serve for asymptotic investigation in these regimes and provide margins as asymptotes; but this may require problem-specific adaptations.
\item Uncertain initial conditions. 
The proposed framework focuses on parametric uncertainties without considering uncertain initial conditions. 
Uncertain initial conditions can contribute significantly to the observation of isolas in practice. 
To tackle this scenario, re-formulating the proposed framework by integrating uncertain states may be considered.
This could be achieved by constructing an augmented boundary-value problem of the local statistical quantity for periodic orbits \cite{Kuehn12,Ahsan24}.
Solutions of the augmented problem can be traced via continuation and used to investigate the impact of uncertain states.
Alternatively, global analysis can be considered by introducing a bounded uncertainty to the state and employing a boundary map approach \cite{Lamb15, Lamb26}.
\end{enumerate}

\dataccess{This article has no additional data.}

\aucontribute{D.H. led the development of the work, with supervisory support from D.A.W.B. and S.A.N.
All authors contributed to the preparation of the manuscript.}

\funding{This work was supported by the EPSRC Programme Grant EP/R006768/1: Digital twins for improved dynamic design.}

\appendix
\section{Parametric Sensitivity of the Initial State}\label{sec:appendixA}
The periodicity~$(\ref{eq:BVE})$ and extremal condition~$(\ref{eq:least_favourable_con_time})$ can be differentiated with respect to uncertainty, yielding
\begin{equation}\label{eq:periodicity_sensitivity}
\frac{{\rm d}\mathbf{x}}{{\rm d} t}\left(T,\mathbf{p}; \mathbf{x}_0\right)\frac{\partial T}{\partial \mathbf{p}}\frac{\partial \mathbf{p}}{\partial \mathbf{\epsilon}} +
\frac{\partial\mathbf{x}}{\partial \mathbf{p}}\left(T,\mathbf{p}; \mathbf{x}_0\right)\frac{\partial \mathbf{p}}{\partial \mathbf{\epsilon}} +
\frac{\partial\mathbf{x}}{\partial \mathbf{x}_0}\left(T,\mathbf{p}; \mathbf{x}_0\right) \frac{\partial\mathbf{x}_0}{\partial \mathbf{p}}\frac{\partial \mathbf{p}}{\partial \mathbf{\epsilon}} -
\frac{\partial\mathbf{x}_0}{\partial \mathbf{p}}\frac{\partial \mathbf{p}}{\partial \mathbf{\epsilon}} =
\mathbf{0} \,,
\end{equation}
\begin{equation}\label{eq:least_fav_sensitivity}
\frac{\partial L_t}{\partial \mathbf{p}}(0, \mathbf{x,p};\mathbf{x}_0)\frac{\partial \mathbf{p}}{\partial \mathbf{\epsilon}}  +
\frac{\partial L_t}{\partial \mathbf{x}}(0, \mathbf{x,p};\mathbf{x}_0)\frac{\partial\mathbf{x}}{\partial \mathbf{x}_0}\left(0, \mathbf{x,p}; \mathbf{x}_0\right)\frac{\partial \mathbf{x}_0}{\partial \mathbf{p}}\frac{\partial \mathbf{p}}{\partial \mathbf{\epsilon}}  =
\mathbf{0} \,.
\end{equation}
Eqs.~$(\ref{eq:periodicity_sensitivity})$ and $(\ref{eq:least_fav_sensitivity})$ can be rearranged in matrix format \cite{Wilkins09},
\begin{align}\label{eq:initial_sensitivity}
\begin{bmatrix}
\mathbf{\mathcal{M} } - \mathbb{I}^{N\times N}				&	\frac{{\rm d}\mathbf{x}}{{\rm d} t}\left(T,\mathbf{p}; \mathbf{x}_0\right) 				\\
\frac{\partial L_t}{\partial \mathbf{x}}(0, \mathbf{x,p};\mathbf{x}_0)	&	0
\end{bmatrix}
\begin{bmatrix}
\mathbf{S}_0				 				\\
\frac{\partial T}{\partial \mathbf{p}}\frac{\partial \mathbf{p}}{\partial \mathbf{\epsilon}}
\end{bmatrix}
=
\begin{bmatrix}
-\mathbf{S}(T, \mathbf{p}; \mathbf{0})				 				\\
-\frac{\partial L_t}{\partial \mathbf{p}}(0, \mathbf{x,p};\mathbf{x}_0)\frac{\partial \mathbf{p}}{\partial \mathbf{\epsilon}}
\end{bmatrix}\,,
\end{align}
where $\mathbf{\mathcal{M}} = \frac{\partial\mathbf{x}}{\partial \mathbf{x}_0}\left(T,\mathbf{p}; \mathbf{x}_0\right)$ is the monodromy matrix of the dynamics, $\frac{\partial T}{\partial \mathbf{p}}\frac{\partial \mathbf{p}}{\partial \mathbf{\epsilon}}$ is the sensitivity of the response period to uncertainty, $\mathbf{S}_0 = \frac{\partial \mathbf{x}_0}{\partial \mathbf{p}}\frac{\partial \mathbf{p}}{\partial \mathbf{\epsilon}}$ is the parametric sensitivity of the initial state, and $\mathbf{S}(T, \mathbf{p}; \mathbf{0})$ is the sensitivity of states to uncertainty at time $t=T$ with a \emph{constant} initial state, i.e.,~\emph{zero} initial state sensitivity with respect to uncertainty.
Solving for $\mathbf{S}_0$ requires $\mathbf{\mathcal{M}}$ and $\mathbf{S}(T, \mathbf{p}; \mathbf{0})$.

The monodromy matrix is governed by the variational dynamics of the system -- the linearised dynamics that describes how small perturbations, $\mathbf{H}$, evolve in the neighbourhood of a given solution, namely,
\begin{equation}\label{eq:first_variant}
\frac{{\rm d}\mathbf{H}}{{\rm d}t} = \mathbf{A}(t, \mathbf{\epsilon}) \mathbf{H}\,,
\end{equation}
where $\mathbf{A} = \frac{\partial\mathbf{F}}{\partial \mathbf{x}}$, defined by Eq.~$(\ref{eq:EoM})$.
The monodromy matrix, $\mathbf{\mathcal{M} } := \mathbf{H}(t=T)$, can be obtained by initialising Eq.~$(\ref{eq:first_variant})$ with $\mathbf{H}(t=0) = \mathbb{I}^{N\times N}$.

To calculate $\mathbf{S}(T, \mathbf{p}; \mathbf{0})$, Eq.~$(\ref{eq:EoM})$ is differentiated with respect to uncertainty to give the sensitivity dynamics
\begin{equation}\label{eq:PDE_par_sensitivity}
\frac{\rm d}{{\rm d} t} \mathbf{S}(t, \mathbf{\epsilon}) = \mathbf{A}(t, \mathbf{\epsilon})\mathbf{S}(t, \mathbf{\epsilon}) + \mathbf{B}(t, \mathbf{\epsilon})\,,
\end{equation}
where $\mathbf{S}= \frac{\partial\mathbf{x}}{\partial\mathbf{p}}\frac{\partial \mathbf{p}}{\partial \mathbf{\epsilon}}\in \mathbb{R}^{N\times M}$ is the parametric sensitivity matrix, $\mathbf{A} = \frac{\partial\mathbf{F}}{\partial \mathbf{x}}$, and $\mathbf{B} = \frac{\partial\mathbf{F}}{\partial \mathbf{p}}\frac{\partial \mathbf{p}}{\partial \mathbf{\epsilon}}$.
 $\mathbf{S}(T, \mathbf{p}; \mathbf{0})$ can be obtained by solving Eq.~$(\ref{eq:PDE_par_sensitivity})$ for $t=T$ with zero initial input.

By substituting $\mathbf{\mathcal{M} }$ (via Eq.~$(\ref{eq:first_variant})$) and $\mathbf{S}(T, \mathbf{p}; \mathbf{0})$ (via Eq.~$(\ref{eq:PDE_par_sensitivity})$) into Eq.~$(\ref{eq:initial_sensitivity})$, the initial state sensitivity, $\mathbf{S}_0$, can be solved.
It is then fed into Eq.~$(\ref{eq:least_favourable_con_uncertainty})$ for determining the extremal conditions.

\section{Rewriting the non-autonomous system}\label{sec:appendixB}
Additional variables are introduced to rewtrite equations of motion~$(\ref{eq:two-mode-EoM})$ in an autonomous state-space form as Eqs.~$(\ref{eq:EoM})$.
Here, variables of a Stuart–Landau oscillator, i.e.,~$x_5$ and $x_6$, are used to replace the periodic excitations.
The resulting state-space description is given by 
\begin{subequations}
\begin{align}
\dot{x}_1 & =   x_2\,, \nonumber \\
\dot{x}_2 & =  \frac{1}{m_1}\left[ - \left(k_1 + k_2 \right)x_1
					   +	k_2 x_3
					   -  \left(c_1 + c_2 \right)x_2
				             +     c_2 x_4
					   - \alpha_1 x_1^3
					  -   \alpha_2 \left(x_1 - x_3 \right)^3
					  + F_1 x_6 \right]\,, \nonumber \\
\dot{x}_3 & =  x_4\,, \nonumber \\
\dot{x}_4 & = \frac{1}{m_2}\left[ - \left(k_2 + k_3 \right)x_3
					   +	k_2 x_1
					   -  \left(c_2 + c_3 \right)x_4
				             +     c_2 x_2
					   - \alpha_3 x_3^3
					  -   \alpha_2 \left(x_3 - x_1 \right)^3
					  + F_2 x_6 \right]\,,  \nonumber \\
\dot{x}_5 & = x_5 +\omega x_6 - x_5 \left(x_5^2 + x_6^2\right)\,, \nonumber\\
\dot{x}_6 & = -\omega x_5 + x_6 - x_6 \left(x_5^2 + x_6^2\right)\,. \nonumber
\end{align}
\end{subequations}

\section{Validation of the results}\label{sec:appendixC}

\begin{figure}[t]
   \centering
   \includegraphics[width=13cm]{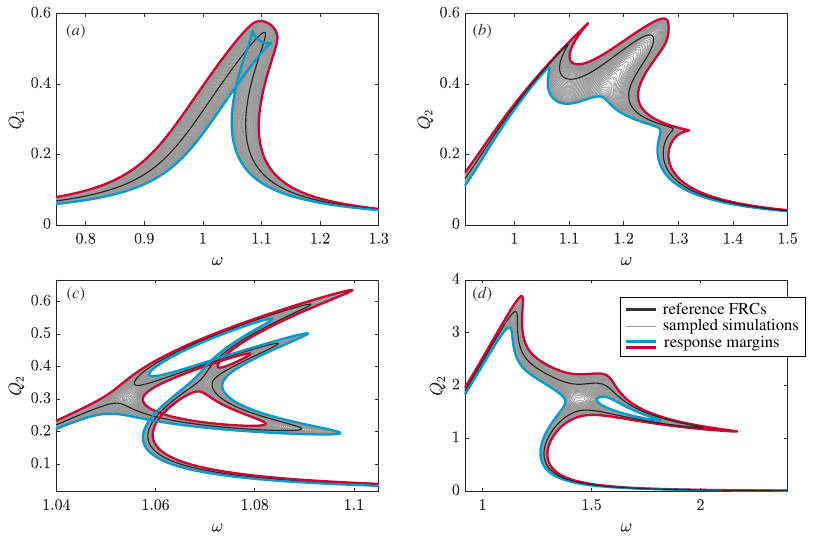}
   \caption{Validation of the UQ results using simulations over a gridded uncertainty domain.
		$(a)$, $(b)$, $(c)$, and $(d)$ correspond to examples in $\S 4$\ref{sec:duffing-like}, \ref{sec:mode-int}, \ref{sec:sym}, and \ref{sec:isola} respectively.
		(Online version in colour.)}
   \label{fig:validation}
\end{figure}

The results in $\S\ref{sec:egs}$ are validated by simulations over the gridded domain of uncertain parameters, shown in Fig.~\ref{fig:validation}, where the simulation results are grey lines whilst the margins, obtained via the proposed method, are blue and red lines.
Panels~$(a)$, $(b)$, $(c)$, and $(d)$ corresponds to examples in $\S 4$\ref{sec:duffing-like}, \ref{sec:mode-int}, \ref{sec:sym}, and \ref{sec:isola} respectively.

\nocite{*}
\bibliographystyle{RS}
\bibliography{bibrefs}

\end{document}